\documentclass[12pt]{article}
\newcommand{\ignore}[1]{}
\pdfoutput=1
\usepackage{amsmath, amsthm, amsfonts, amssymb, mathrsfs, eufrak, pgf}

\usepackage{graphicx}

\usepackage{float}
\newfloat{figure}{H}{lof}
\floatname{figure}{\figurename}

\numberwithin{figure}{section}

\DeclareMathAlphabet{\eufrak}{U}{}{}{}  
\SetMathAlphabet\eufrak{normal}{U}{euf}{m}{n}
\SetMathAlphabet\eufrak{bold}{U}{euf}{b}{n}

\newenvironment{Proof}{\removelastskip\par\medskip
\noindent{\em Proof.}
\rm}{\hfill$\square$\\\par\medbreak}

\def\SURE{{\mathord{{\rm SURE \ \!}}}}

\def\real{{\mathord{{\rm I\kern-2.8pt R}}}}        
\def\inte{{\mathord{{\rm I\kern-2.8pt N}}}}
\def\PP{{\mathord{{\rm I\kern-2.8pt P}}}}

\def\real{{\mathord{\mathbb R}}}

\def\inte{{\mathord{\mathbb N}}}

\def\Cov{{\mathrm{{\rm Cov}}}}

\def\crb{{\mathrm{{\rm R}}}}

\newcommand{\sgn}{\mathrm{sign }}

\def\P{\mathbb{P}}
\def\E{\mathop{\hbox{\rm I\kern-0.20em E}}\nolimits}

\def\Cov{\mathop{\hbox{\rm Cov}}\nolimits}

\newtheorem{prop}{Proposition}[section]

\newtheorem{lemma}[prop]{Lemma}
\newtheorem{definition}[prop]{Definition}

\def\Dom{{\mathrm{{\rm Dom \ \! }}}}

\def\var{{\mathrm{{\rm Var \ \! }}}}
\def\argmin{{\mathrm{{\rm argmin}}}}

\numberwithin{equation}{section}

\title{
\Huge 
 SURE shrinkage of Gaussian paths and signal identification 
} 

\author
{
Nicolas Privault\footnote{nprivaul@cityu.edu.hk}
\\
Department of Mathematics
\\
City University of Hong Kong
\\
Tat Chee Avenue
\\
Kowloon Tong
\\
Hong Kong
\and
Anthony R\'eveillac\footnote{anthony.reveillac@univ-lr.fr}
\\
Institut f\"ur Mathematik
\\
Humboldt-Universit\"at zu Berlin
\\
Unter den Linden 6
\\
10099 Berlin
\\
Germany
}

\allowdisplaybreaks

\begin{document}
\hyphenation{func-tio-nals}

\maketitle
 
\vspace{-0.4cm} 

\begin{abstract} 
 Using integration by parts on Gaussian space 
 we construct a Stein Unbiased Risk Estimator (SURE) 
 for the drift of Gaussian processes using their 
 local and occupation times. 
 By almost-sure minimization of the SURE risk of 
 shrinkage estimators we derive an estimation and de-noising 
 procedure for an input signal perturbed by a 
 continuous-time Gaussian noise. 
\end{abstract} 

\normalsize

\vspace{0.5cm}


\normalsize

\baselineskip0.7cm

\section{Introduction}
\noindent 
 Let $X$ be a Gaussian random vector on $\real^d$ with unknown mean $m$ 
 and known covariance matrix $\sigma^2 \textbf{I}_d$ 
 under a probability measure $\P_m$. 

 It is well-known \cite{stein} that given 
 $g:\real^d\to\real^d$ a sufficiently smooth function, 
 the mean square risk $\|X+g(X)-m \|_{\real^d}^2$ 
 of $X+g(X)$ to $m$ can be estimated unbiasedly by 
\begin{equation} 
\label{eq:SteinRisk} 
 \SURE : = 
 \sigma^2 
 d + \sum_{i=1}^d g_i(X)^2 + 2 \sum_{i=1}^d \nabla_i g(X)
, 
\end{equation} 
 from the identity 
\begin{equation} 
\label{l1} 
 \E_m \left[\|X+g(X)- m \|_{\real^d}^2\right] = \sigma^2 
 d 
 + 
 \E_m \left[ 
 \sum_{i=1}^d g_i(X)^2 
 + 
 2 
 \sum_{i=1}^d 
 \nabla_i g(X)\right]
\end{equation}  
 which is obtained by Gaussian integration by parts 
 under $\P_m$. 
 The estimator \eqref{eq:SteinRisk}, 
 which is independent of $m$, 
 is called the Stein Unbiased Risk Estimate (SURE). 
\\ 
 
 When $(g^\lambda)_{\lambda \in \Lambda}$ is a 
 family of functions it makes sense to almost surely minimize 
 the {\em Stein Unbiased Risk Estimate} \eqref{eq:SteinRisk} 
 of $g^\lambda$ with respect to the parameter $\lambda$. 
 This point of view has been developed by Donoho and 
 Johnstone~\cite{DonohoJohnstone} for the design of 
 spatially adaptive estimators 
 by shrinkage of wavelet coefficients of noisy data 
 via 
$$ 
 X + g^\lambda ( X ) = \lambda \eta ( X / \lambda )
, 
$$ 
 where $\eta (x)$ is a threshold function. 
\\ 
 
 In this paper we construct a Stein type Unbiased Risk Estimator
 for the deterministic drift $(u_t)_{t\in \real_+}$ 
 of a one dimensional Gaussian processes $(X_t)_{t\in [0,T]}$ 
 via an extension of the identity \eqref{l1} introduced in 
 \cite{pr3}, \cite{p-r-c} on the Wiener space. 
 For example, given $\alpha (t)$ and $\lambda (t)$ two functions given in 
 parametric form, the SURE risk of the estimator 
$$ 
 X_t
 + 
 \xi^{\alpha , \lambda}_t (X_t) 
 = 
 \alpha (t) + 
 \lambda (t) 
 \eta_S \left( 
 \frac{X_t-\alpha (t)}{ \lambda (t) } 
 \right) 
, 
 \qquad 
 t\in [0,T], 
$$ 
 where $\eta_H$ is the hard threshold function 
 \eqref{hthre} below, is given by 
$$ 
 \SURE (X+\xi^{\alpha , \lambda }(X) ) 
 = 
 T 
 + 
 \int_0^T 
 \frac{( X_t -  \alpha (t) )^2}{\gamma (t,t)} 
 {\bf 1}_{\{ \vert X_t - \alpha (t) 
 \vert \leq \lambda \sqrt{\gamma(t,t)} \} }
 dt 
 + 
 2 
 \lambda 
 \bar{\ell}_T^\lambda 
 - 
 2 
 \bar{L}_T^\lambda 
, 
$$ 
 where $\gamma (s,t) = \Cov ( X_s , X_t )$, 
 $0\leq s , t \leq T$, denotes the covariance of $(X_t)_{t\in [0,T]}$ 
 and $\bar{\ell}_T^\lambda$, $\bar{L}_T^\lambda$ 
 respectively denote the local and occupation time of 
$$ 
 ( | X_t - \alpha (t) |/\sqrt{\gamma (t,t)})_{t\in [0,T]} 
, 
$$ 
 cf. Proposition~\ref{Prop:EqualityOfRisks.2}. 
 We apply this technique to de-noising and identification of 
 the input signal in a Gaussian channel 
 via the minimization of $\SURE (X+\xi^{\alpha , \lambda }(X) )$. 
 This yields in particular an estimator 
 of the drift of $X_t$ from the estimation of 
 $\alpha (t)$, and 
 an optimal noise removal threshold from the estimation 
 of $\lambda$. 
 This approach differs from classical signal detection techniques 
 which usually rely 
 on likelihood ratio tests, cf e.g. \cite{poor}, Chapter~VI. 
 It also requires an a priori hypothesis on the parametric 
 form of $\alpha (t)$. 
\\ 
 
 We proceed as follows. 
 In Section~\ref{1.1.1} we recall our framework 
 of functional estimation of drift trajectories. 
 In Section~\ref{3} we derive Stein's unbiased risk estimate 
 for the estimation of the drift of Gaussian processes. 
 In Section~\ref{4} we discuss its application to 
 soft thresholding for Gaussian 
 processes using the local time and obtain an upper 
 bound for the risk of such estimators. 
 We also show the existence of an optimal parameter and 
 the smoothness of the risk function. 
 In Section~\ref{4.1} we consider the case of 
 hard thresholding. 
 In Section~\ref{5} we consider several numerical examples 
 where $\alpha (t)$ is given in parametric form. 
 In Section~\ref{2} we recall some elements of 
 stochastic analysis of Gaussian processes. 
\section{Functional drift estimation} 
\label{1.1.1} 
 In this section we recall the setting of functional 
 drift estimation to be used in this paper. 
 Given $T>0$ we 
 consider a real-valued centered Gaussian process $X = (X_t)_{t\in [0,T]}$
 with non-vanishing covariance function 
$$ 
 \gamma (s,t) = \E [ X_s X_t ], \qquad s,t\in [0,T],
$$ 
 on a probability space $(\Omega , {\cal F} , \P )$,
 where $({\cal F})_{t\in [0,T]}$ 
 is the filtration generated by $(X_t)_{t\in [0,T]}$. 
 Assume that under a probability measure $\P_u$ 
 we observe the paths of $(X_t)_{t\in [0,T]}$ decomposed as
$$ 
 X_t=u_t+X_t^u, \quad t \in [0,T], 
$$ 
 where 
 $u=(u_t)_{t\in [0,T]}$ is a square integrable $\mathcal{F}$-adapted process and $(X_t^u)_{t \in [0,T]}$ is a centered Gaussian process with covariance 
$$ 
 \gamma (s,t) = \E_u [ X^u_s , X^u_t ], \qquad 0 \leq s , t \leq T, 
$$ 
 where $\E_u$ denotes the expectation under $\P_u$. 
 Given a continuous time observation of the process $(X_t)_{t\in [0,T]}$ 
 we will propose estimators of the unknown drift function $u$. 
\begin{definition}
The risk of an estimator $\xi:=(\xi_t)_{t\in [0,T]}$ to $u$ is defined as
$$
 R ( \gamma , \mu , \xi ) : = 
 \E_u \left[
 \int_0^T  
 | \xi_t - u_t |^2 \mu (dt) 
 \right] 
$$ 
 where $\mu$ is a positive measure on $[0,T]$.  
\end{definition}
 Examples of risk measures $\mu$ include the Lebesgue measure and 
\begin{equation} 
\label{m0} 
 \mu(dt) 
 = 
 \sum_{i=1}^n 
 a_i 
 \delta_{t_i}(dt) 
, 
 \qquad 
 a_1,\ldots ,a_n > 0 
, 
\end{equation} 
 in which case the risk of the estimator is computed 
 from discrete values of the sample path observed at times 
 $t_1,\ldots,t_n$, $n\geq 1$. 
\begin{definition}
A drift estimator $(\xi_t)_{t\in [0,T]}$ is called unbiased if 
$$
 \E_u [\xi_t]= \E_u [u_t], \quad t \in [0,T],
$$
 for all square-integrable ${\cal F}_t$-adapted process
 $(u_t)_{t\in [0,T]}$, where 
 $({\cal F}_t)_{t\in [0,T]}$ is the filtration generated
 by $(X_t)_{t\in [0,T]}$.
\end{definition}
\noindent
 In the sequel we will consider the canonical process $(X_t)_{t\in [0,T]}$ 
 as an unbiased estimator $\hat{u} : = (X_t)_{t\in [0,T]}$ 
 of its own drift $(u_t)_{t\in [0,T]}$
 under $\P_u$, with risk 
$$ 
 \crb (\gamma ,\mu, \hat{u} ) 
 : = 
 \E_u \left[ \int_0^T | X_t - u_t|^2 \mu ( dt ) \right] 
 = \int_0^T \gamma( t,t) \mu ( dt ) 
$$ 
 Recall that the estimator $\hat{u} = (X_t)_{t\in [0,T]}$ is minimax i.e. 
$$ 
 \crb (\gamma ,\mu, \hat{u} ) 
 = 
 \inf_\xi 
 \sup_{v \in \Omega } 
 \E_v \left[ \int_0^T | \xi_t - v_t|^2 \mu ( dt ) \right] 
, 
$$ 
 cf. Proposition 3.2 of \cite{pr3}. 
 In addition, when $(X_t)_{t\in [0,T]}$ has independent increments 
 and $(u_t)_{t\in [0,T]} \in L^2 (\Omega \times [0,T] , \P_u \otimes \mu )$ 
 is square-integrable and adapted, then for any adapted and unbiased estimator $\xi$ the Cramer-Rao bound 
\begin{equation} 
\label{11} 
 \E_u \left[ \int_0^T |\xi_t - u_t |^2  \mu ( dt ) \right] 
 \geq 
 \crb (\gamma,\mu, \hat{u} ) 
, 
\end{equation} 
 holds for any unbiased and adapted estimator $(\xi)_{t\in [0,t]}$ of 
 $(u(t))_{t\in [0,T]} \in L^2 (\Omega \times [0,T] , \P_u \otimes \mu )$  
 and is attained by $\hat{u}$, cf. Proposition~4.3 of \cite{pr3}, 
 hence $\hat{u} = (X_t)_{t\in [0,T]}$ 
 is an efficient estimator of its own drift $u$. 
\section{Stein's unbiased risk estimate} 
\label{3} 
 Instead of using the minimax estimator $\hat{u}$ 
 we will estimate the drift of $(X_t)_{t\in [0,T]}$ 
 by the almost sure minimization of a 
 Stein Unbiased Risk Estimator for Gaussian processes, 
 constructed in the next proposition by analogy with \eqref{eq:SteinRisk}. 
 In the next proposition we use the gradient operator $D_t$ 
 whose definition and properties are recalled in the 
 appendix, cf. Definition~\ref{defd} and Lemma~\ref{prt}. 
\begin{prop} 
\label{lemma1} 
 For any $(\xi_t)_{t\in [0,T]} \in 
 L^2 (\Omega \times [0,T] , \P_u \otimes \mu )$ such that
 $\xi_t \in \Dom (\nabla)$, $t\in [0,T]$, and
 $(D_t\xi_t)_{t\in [0,T]}\in L^1 (\Omega \times [0,T] , \P_u \otimes \mu )$, 
 the quantity 
\begin{equation}
\label{eq:defiSURE}
 \SURE_\mu (X+\xi):= 
 \crb (\gamma ,\mu, \hat{u} )+\Vert \xi \Vert_{L^2 ( [0,T] , d \mu )}^2+ 2 
 \int_0^T
 D_t \xi_t
 \mu ( dt ) 
\end{equation}
 is an unbiased estimator of the mean square risk 
 $\Vert X + \xi - u \Vert_{L^2([0,T] , d\mu )}^2$. 
\end{prop}
\begin{Proof}
 From Lemma~\ref{prt} we have 
\begin{eqnarray*}
\nonumber 
\lefteqn{
 \E_u \left[
 \Vert X + \xi - u \Vert_{L^2([0,T] , d\mu )}^2 \right]
 =
 { \E_u \left[ \int_0^T \Big{|} X^u_t + \xi_t \Big{|}^2
 \mu ( dt )
 \right] }
}
\\
\nonumber 
 &=& { \E_u \left[ \int_0^T | X^u_t |^2
 \mu ( dt )
 \right] +
 \E_u \left[ 
 \Vert \xi \Vert_{L^2 ( [0,T] , d \mu )}^2 
 \right] 
 +
 2 \E_u \left[ \int_0^T X^u_t \xi_t
 \mu ( dt )
 \right] } 
\\
\nonumber 
 &=&
 \crb (\gamma ,\mu, \hat{u} )
+
 \E_u \left[ 
 \Vert \xi \Vert_{L^2 ( [0,T] , d \mu )}^2 
 \right] 
 + 2 \E_u \left[ \int_0^T D_t\xi_t
 \mu ( dt )
 \right]
\\
 &=&
 \E_u 
 \left[ 
 \SURE_\mu (X+\xi) 
 \right] 
.
\end{eqnarray*}
\end{Proof} 
 Unlike the pointwise mean square risk 
 $\Vert X + \xi - u \Vert_{L^2([0,T] , d\mu )}^2$, 
 the SURE risk estimator does not depend on the estimated 
 parameter $u$. 
\\ 
 
 Given a family $(\xi^{\lambda})_{\lambda \in \Lambda}$ 
 of estimators indexed by a parameter space $\Lambda$, 
 we consider the estimator $X+\xi^{\lambda^*}$ that 
 almost-surely minimizes the SURE risk, with 
$$ 
 \lambda^* 
 = 
 \argmin_{\lambda \in \Lambda} 
 \SURE_\mu (X+\xi^\lambda ) 
. 
$$ 
 For all values of $\lambda$ 
 the SURE risk of the estimator $X+\xi^{\lambda^*}$ improves on the 
 mean square risk of $X+\xi^{\lambda}$. 
 
 Precisely for all $\nu \in \Lambda$ we have 
\begin{eqnarray*} 
 \E_u [ \SURE_\mu (X+\xi^{\lambda^*} ) ] 
 & \leq & 
 \E_u [ \SURE_\mu (X+\xi^\nu ) ] 
\\ 
 & = & 
 \E_u \left[
 \Vert \xi^\nu - u \Vert^2_{L^2([0,T], \mu)} 
 \right] 
\\ 
 & = & 
 \inf_{\lambda} 
 \E_u \left[
 \Vert \xi^\lambda - u \Vert^2_{L^2([0,T], \mu)} 
 \right] 
. 
\end{eqnarray*} 

 In the sequel we will apply the above to a process 
 $(\xi_t)_{t\in [0,T]}$ given as a funtion 
 $\xi_t = \xi_t (X_t)$ of $X_t$, $t\in [0,T]$. 
 In particular we will discuss estimation and thresholding 
 for estimators of the form 
\begin{equation} 
\label{1} 
 X_t 
 + 
 \xi^{\alpha , \lambda}_t (X_t) 
 = 
 \alpha (t) + \lambda (t) 
 \eta \left( 
 \frac{X_t-\alpha (t)}{\lambda (t)} 
 \right) 
, 
\end{equation} 
 where $\eta : \real \to \real$ is a threshold function 
 with support in $(-\infty , - 1 ] \cup [ 1 , \infty )$. 
\\ 
 
 In particular we will apply our method to the joint estimation 
 of parameters $\alpha$, $\lambda$, successively 
 in case 
 $\alpha (t) = \alpha$, 
 $\alpha (t) = \alpha t$, 
 and 
 $\lambda (t) = \lambda \sqrt{\gamma (t,t)}$. 
\section{Soft threshold} 
\label{4} 
 In this section we construct an example of 
 SURE shrinkage by soft thresholding 
 in the framework of Proposition~\ref{lemma1}, 
 with application to identification and 
 de-noising in a Gaussian signal. 
 In case $\eta$ is the soft threshold function 
\begin{equation} 
\label{sthre} 
 \eta_S (y)= \sgn ( y ) ( | y | -1)^+, \quad y \in \real 
, 
\end{equation} 
 the function $\xi^{\alpha , \lambda}_t$ in \eqref{1} 
 becomes 
$$ 
 \xi^{\alpha , \lambda}_t ( x ) 
 = 
 - \sgn( x - \alpha (t) ) \min ( 
 \lambda (t) , \vert x - \alpha (t) \vert ) 
, \qquad x \in \real 
, 
$$ 
 where $\lambda (t) \geq 0$ is a given level function. 
\begin{prop} 
\label{Prop:EqualityOfRisks} 
 We have $\P$-a.s 
\begin{eqnarray} 
\label{eq:EqualityOfRisks} 
\lefteqn{ 
 \SURE_\mu (X+\xi^{\alpha , \lambda } (X))} 
\\ 
\nonumber 
 & = & 
 \crb (\gamma ,\mu, \hat{u} ) +  \int_0^T 
 \vert X_t - \alpha (t) \vert^2 \wedge \lambda^2 (t) 
 \mu(dt)-2\int_0^T 
 {\bf 1}_{\{ \vert X_t - \alpha (t) \vert \leq \lambda (t) \} 
 } \gamma(t,t) \mu(dt) 
. 
\end{eqnarray}
\end{prop} 
\begin{Proof} 
 Since 
 $\frac{d}{dx} \xi^{\alpha , \lambda}_t ( x ) = 
 - {\bf 1}_{\{ \vert x - \alpha (t) \vert \leq \lambda (t) \}}$, 
 we have 
\begin{eqnarray*}
 \int_0^T 
 D_t \xi^{\alpha , \lambda}_t
 ( 
 X_t
 ) 
 \mu ( dt )
 &=& 
 - 
 \int_0^T 
 {\bf 1}_{\{ \vert X_t - \alpha (t) \vert \leq \lambda (t) \} 
 } 
 D_t 
 X_t 
\mu(dt)\nonumber\\
&=& 
 - 
 \int_0^T 
 {\bf 1}_{\{ \vert X_t - \alpha (t) \vert \leq \lambda (t) \} 
 } 
 \gamma(t,t) 
 \mu(dt), 
\end{eqnarray*}
 hence the conclusion from Proposition~\ref{lemma1}. 
\end{Proof}
 The risk associated to discrete observations $(X_{t_1},\ldots,X_{t_n})$ can 
 be computed via Proposition~\ref{Prop:EqualityOfRisks} by choosing the 
 risk measure \eqref{m0}, in which case 
 Relation~\eqref{eq:EqualityOfRisks} becomes 
\begin{eqnarray*}
\lefteqn{ 
 \SURE ( 
 X+\xi^{\alpha , \lambda} (X) 
 ) 
} 
 \\ 
 & = & 
 \crb (\gamma ,\mu, \hat{u} )
 +\sum_{i=1}^n 
 \vert X_{t_i} - \alpha (t_i) \vert^2 \wedge \lambda^2 (t_i) 
 - 
 2 \sum_{i=1}^n \gamma(t_i,t_i) 
 {\bf 1}_{ \{ 
 \vert X_{t_i} - \alpha (t_i) \vert \leq \lambda (t_i) 
 \} 
 } 
. 
\end{eqnarray*} 
 which is analog to the finite dimensional SURE risk 
\begin{equation} 
\label{eq:SURERisk}
 \SURE ( X + g^\lambda (X) 
 ) 
 = 
 d +\sum_{i=1}^d (\vert x_i\vert \wedge \lambda )^2-2 \#\{i; \; \vert x_i\vert \leq \lambda \} 
\end{equation} 
 of \cite{DonohoJohnstone2}. 
 In the simulations of Section~\ref{5} 
 we effectively use such risk measures when discretizing 
 the signal. 
 More precisely, when $\mu (dt) = f(t) dt$ has a density 
 $f(t)$ with respect to the Lebesgue measure and 
 $$ 
 \mu_n(dt)=\sum_{i=1}^{n-1} 
 f(t_i) (t_{i+1}-t_i) \delta_{t_i}(dt)
, 
$$ 
 Relation~\eqref{eq:EqualityOfRisks} shows that 
 $\SURE_{\mu_n} (X+\xi^{\alpha,\lambda}(X))$ 
 becomes a consistent estimator of the risk 
 $\SURE_\mu (X+\xi^{\alpha,\lambda}(X))$ as $n$ goes to infinity.
\\ 

\noindent 
 Taking 
$$ 
 \mu(dt) = \gamma^{-1} (t,t) dt 
 \quad
 \mbox{and} 
 \quad 
 \lambda (t) = \lambda \sqrt{\gamma (t,t)}, 
 \qquad 
 \lambda >0, \quad 
 t\in [0,T], 
$$  
 and letting 
\begin{equation} 
\label{lett} 
 \bar{L}_T^\lambda 
 : = 
 \int_0^T {\bf 1}_{ \{ 
 \vert X_t - \alpha (t) \vert \leq \lambda \sqrt{\gamma(t,t)} 
 \} 
 } dt 
\end{equation} 
 denote the occupation time of the process 
$$ 
 Z_t^{\alpha,\gamma} : = \frac{X_t-\alpha (t) }{\sqrt{\gamma (t,t)}}, 
 \qquad t\in [0,T], 
$$ 
 up to time $T$ in the set $[ - \lambda , \lambda ]$, 
 Proposition~\ref{Prop:EqualityOfRisks} yields the 
 identity 
\begin{equation} 
\label{1.1} 
 \SURE_\mu ( X+\xi^{\alpha , \lambda} (X) ) 
 = 
 T 
 + 
 \int_0^T \left( 
 | 
 Z_t^{\alpha,\gamma} 
 | 
 \wedge \lambda\right)^2 dt 
 - 
 2 
 \bar{L}_T^\lambda 
. 
\end{equation}  
 As a consequence we obtain the following bound for the 
 risk of the thresholding estimator $X+\xi^{\alpha , \lambda} (X)$. 
\begin{prop} 
 Assume that $u \in L^2([0,T] , d \mu )$ is a deterministic 
 function and let $\mu(dt):=\gamma(t,t)^{-1} dt$. 
 Then for all fixed $\lambda \geq 0$ we have 
$$ 
 \E_u [ 
 \Vert X + \xi^{\alpha , \lambda} (X) - u \Vert_{L^2([0,T] , d\mu )}^2 
 ] 
 \leq 
 ( 1 + \lambda^2 ) 
 \left(T\wedge \int_0^T  
 | 
 u(t) - \alpha (t) 
 |^2  
 \mu ( dt ) \right)
 + 
 T 
 ( 
 1 + \lambda 
 ) 
 e^{-\frac{\lambda^2}{2}} 
. 
$$
\end{prop} 
\begin{Proof} 
 We have 
$$ 
 \SURE_\mu ( X+\xi^{\alpha , \lambda} (X) ) 
 = 
 T 
 + 
 \int_0^T \left( 
 | 
 Z_t^{\alpha,\gamma} 
 | 
 \wedge \lambda\right)^2 dt 
 - 
 2 
 \int_0^T {\bf 1}_{ \{ 
 \vert X_t - \alpha (t) \vert \leq \lambda \sqrt{\gamma(t,t)} 
 \} 
 } dt 
$$ 
 hence 
$$ 
 \E_u [ 
 \SURE_\mu ( X+\xi^{\alpha , \lambda} (X) ) 
 ] 
 \leq T ( 1 + \lambda^2)
, 
$$
 and 
\begin{eqnarray*} 
 \E_u [ 
 \SURE_\mu ( X+\xi^{\alpha , \lambda} (X) ) 
 ] 
 & \leq & 
 \int_0^T 
 1+\E_u [ 
 | 
 Z_t^{\alpha,\gamma} 
 |^2]\wedge \lambda^2  
 - 
 2  
 \P_u 
 ( 
 \vert 
 Z_t^{\alpha,\gamma} 
 \vert 
 \leq \lambda 
 ) 
 dt 
\\ 
 & \leq & 
 \int_0^T 
 (1+\lambda^2) \left(e^{-\frac{\lambda^2}{2}}+\frac{\vert u(t)-\alpha(t) \vert^2}{\gamma(t,t)}\right) 
 dt 
\\ 
 & \leq & 
 ( 1 + \lambda^2 ) 
 \int_0^T 
 | 
 u(t) - \alpha (t) 
 |^2 
 \mu ( dt ) 
 + 
 T 
 ( 
 1 + \lambda^2 
 ) 
 e^{-\frac{\lambda^2}{2}} 
,
\end{eqnarray*} 
 where we recall that from \cite{DonohoJohnstone2}, Appendix 1, 
 we have for every $t$ in $[0,T]$ that
$$ 1+\E_u [\vert Z_t^{\alpha,\gamma} \vert^2]\wedge \lambda^2-2\P_u(\vert Z_t^{\alpha,\gamma} \vert \leq \lambda) \leq (1+\lambda^2) \left( e^{-\frac{\lambda^2}{2}} + \frac{\vert u(t)-\alpha(t)\vert^2}{\gamma(t,t)}\right)$$
and we conclude from Proposition~\ref{lemma1}. 
\end{Proof} 
 From this proposition it follows that 
 $\SURE_\mu ( X+\xi^{\alpha , \lambda} (X) )$ 
 is independent of large values 
 $\Vert u - \alpha \Vert_{L^2([0,T])}$, 
 while its growth at most as $1+\lambda^2$ 
 in $\lambda \geq 0$. 
\\ 

 Since $\lambda \mapsto \SURE_\mu ( X+\xi^{\alpha ,\lambda} (X) )$ in 
 \eqref{1.1} is lower bounded 
 by $-T$ and equal to $0$ when $\lambda = 0$, the optimal threshold 
\begin{equation} 
\label{argm} 
 \lambda^* := \argmin_{\lambda} \SURE_\mu ( X+\xi^{\alpha , \lambda} (X) ) 
\end{equation} 
 exists almost surely in $[0,\infty )$. 
\\ 

 In addition we have the following proposition which important for the 
 numerical search of an optimal parameter value. 
\begin{prop} 
\label{minnum} 
 The function $\lambda \mapsto 
 \SURE_\mu ( X + \xi^{\alpha , \lambda} (X) )$ 
 is continuously differentiable. 
\end{prop} 
\begin{Proof} 
 Letting 
$$ 
 \Delta (s,t) = \var_u ( Z^{\alpha,\gamma}_t - Z^{\alpha,\gamma}_s) = 
 2 - 2 \frac{\gamma (s,t)}{\sqrt{\gamma (s,s)\gamma (t,t)}}, 
 \qquad 
 0\leq s , t \leq T 
, 
$$ 
 under Condition~\eqref{eq:TL}, the local time 
$$ 
 \bar{\ell}_T^\lambda : = 
 \frac{d}{d\lambda} \bar{L}_T^\lambda 
$$ 
 of $(|Z^{\alpha ,\gamma}_t|)_{t\in [0,T]}$ exists almost surely,  
 cf. Section~\ref{2}, and we have 
\begin{eqnarray*} 
\frac{\partial }{\partial \lambda} 
 \SURE_\mu ( X + \xi^{\alpha , \lambda} (X) ) 
 & = & 
 \frac{\partial }{\partial \lambda} 
 \int_0^T 
 \left( 
 | 
 Z^{\alpha , \gamma}_t 
 | 
 \wedge \lambda\right)^2 dt-2 \bar{\ell}_T^\lambda 
\\
&=& 
 2\lambda \int_0^T {\bf 1}_{ \{  
 \vert X_t - \alpha (t) \vert \geq \lambda \sqrt{\gamma(t,t)} 
 \} 
 } dt-2\bar{\ell}_T^\lambda 
\\
 &=& 
 2\lambda (T-\bar{L}_T^\lambda ) -2\bar{\ell}_T^\lambda, 
\end{eqnarray*}
 which is a continuous function of $\lambda$ since the
 covariance $\gamma (s,t)$ does not vanish, cf. e.g. 
 Theorem~26.1 of \cite{GemanHorowitz}. 
\end{Proof} 
 Consequently we have 
$$ 
 \frac{\partial}{\partial \lambda} 
 \SURE_\mu ( X + \xi^{\alpha , \lambda } (X) 
 )_{\vert \lambda=0}=-2 \bar{\ell}_T^0, 
$$ 
 hence $\lambda^* > 0$ a.s. when $\ell_T^0$ 
 is a.s. positive, which is the case for example 
 when $X_t$ is a Brownian motion, see Corollary~2.2 of 
 page 240 of \cite{revuz3}, Chapter~VI. 
\\ 

 In practice we will compute $\lambda^*$ numerically 
 by minimization of $\lambda \mapsto 
 \SURE_\mu ( X + \xi^{\alpha ,\lambda} (X) )$ over 
 $\lambda$ in a range $\Lambda = [0,C(T)]$ where $C(T)$ 
 is such that 
$$ 
 \lim_{T \to \infty} 
 \P_u 
 \left( 
 \sup_{t\in[0,T]} | Z^{\alpha,\gamma}_t | \leq C(T) 
 \right) 
 = 
 1. 
$$ 
 This condition is analog to Condition~(31) in \cite{DonohoJohnstone2} 
 and allows us to restrict the range of $\lambda$ when searching 
 for an optimal threshold. 
\\ 

 The function $\alpha (t)$ can be given in parametric 
 form, in which case the parameters will 
 be used to minimize $\SURE_\mu ( X + \xi^{\alpha , \lambda} (X) )$, 
 cf. Section~\ref{5}. 
\section{Hard threshold} 
\label{4.1} 
 Here we use the threshold function 
\begin{equation} 
\label{hthre} 
 \eta_H (y)= y {\bf 1}_{\{ | y | > 1 \} }, 
 \quad y \in \real 
, 
\end{equation} 
 hence 
$$ 
 \xi^{\alpha , \lambda}_t ( x ) 
 = 
 - ( 
 x 
 - 
 \alpha (t) 
 ) 
 {\bf 1}_{\{ 
 \vert x - \alpha (t) \vert < \lambda \sqrt{\gamma(t,t)}  
 \} }, 
 \qquad x \in \real 
, 
$$ 
 where $\lambda \geq 0$ is a level parameter. 
\\ 
 
 In finite dimensions \cite{DonohoJohnstone2} 
 the SURE estimator \eqref{eq:SteinRisk} can not be computed 
 due to the non-differentiability of $\eta_H$, 
 however a deterministic optimal threshold equal to $\sqrt{2\log d}$ 
 can be obtained by other methods, cf. Theorem~4 of 
 \cite{DonohoJohnstone2}. 
\\ 
 
 In continuous time the situation is different 
 due to the smoothing effect of the integral over time. 
 In the next proposition we compute the SURE risk 
 using the local time of Gaussian 
 processes when $\mu(dt) = \gamma^{-1} (t,t) dt$. 
\begin{prop} 
\label{Prop:EqualityOfRisks.2} 
 We have $\P$-a.s 
\begin{eqnarray} 
\nonumber 
 \SURE_\mu (X+\xi^{\alpha , \lambda }(X) ) 
 & = & 
 T 
 + 
 \int_0^T 
 \frac{( X_t -  \alpha (t) )^2}{\gamma (t,t)} 
 {\bf 1}_{\{ \vert X_t - \alpha (t) 
 \vert \leq \lambda \sqrt{\gamma(t,t)} \} }
 dt 
 + 
 2 
 \lambda 
 \bar{\ell}_T^\lambda 
 - 
 2 
 \bar{L}_T^\lambda 
. 
\\ 
& & 
\label{pl} 
\end{eqnarray}
\end{prop}
\begin{proof} 
 Let $\phi\in {\cal C}^\infty_c ([-1,1])$, 
 $\phi \geq 0$ be symmetric around the origin, 
 such that 
$\int_{-1}^1 \phi ( x ) d x = 1$, and let 
$$\phi_\varepsilon ( x ) = \varepsilon^{-1} \phi (\varepsilon^{-1} x ), 
 \qquad  x \in \real, 
 \quad \varepsilon > 0. 
$$ 
 Let 
$$ 
 \xi^{\alpha,\lambda,\varepsilon}_t (x) 
 = 
 \phi_{\varepsilon \sqrt{\gamma (t,t)}} * \xi_t^{\alpha,\lambda} (x) 
 = \int_{-\infty}^\infty \phi_{\varepsilon \sqrt{\gamma (t,t)}} 
 ( y ) \xi_t^{\alpha,\lambda} (x- y ) 
 d y 
, 
$$ 
 denote the convolution of 
 $\phi_{\varepsilon \sqrt{\gamma (t,t)}}$ with $\xi^{\alpha , \lambda}_t$, 
 with 
\begin{eqnarray*} 
 \frac{d}{dx} 
 \phi_{\varepsilon \sqrt{\gamma (t,t)}} 
 * \xi_t^{\alpha,\lambda} (x) 
 & = & 
 \phi_{\varepsilon \sqrt{\gamma (t,t)}} * 
 \frac{d}{dx} 
 {\xi_t^{\alpha,\lambda}} (x) 
\\ 
 & = & 
 \lambda \sqrt{\gamma (t,t)} 
 \phi_{\varepsilon \sqrt{\gamma (t,t)}} 
 ( - \lambda \sqrt{\gamma (t,t)} 
 + x - \alpha (t) 
 ) 
\\ 
 & & 
 + 
 \lambda \sqrt{\gamma (t,t)} 
 \phi_{\varepsilon \sqrt{\gamma (t,t)}} 
 ( \lambda \sqrt{\gamma (t,t)} ) 
 + x - \alpha (t) 
\\ 
 & & 
 - 
 \int_{-\infty}^\infty 
 \phi_{\varepsilon \sqrt{\gamma (t,t)}} 
 ( y ) 
 {\bf 1}_{\{ 
 \vert 
 x - y - \alpha (t) 
 \vert 
 <  
 \lambda \sqrt{\gamma (t,t)} 
 \} }
 dy 
. 
\end{eqnarray*} 
 From the occupation time density formula \eqref{otdf} we have 
\begin{eqnarray} 
\label{eq:RiskXiBis}
\nonumber
 \int_0^T
 D_t \xi^{\alpha , \lambda , \varepsilon}_t (X_t) 
 \mu ( dt )
 &=& 
 \lambda 
 \int_0^T 
 \sqrt{\gamma (t,t)} 
 \phi_{\varepsilon \sqrt{\gamma (t,t)}} 
 ( - \lambda \sqrt{\gamma (t,t)} + X_t - \alpha (t) ) 
 dt
\\ 
\nonumber 
 & & 
 + 
 \lambda 
 \int_0^T 
 \sqrt{\gamma (t,t)} 
 \phi_{\varepsilon \sqrt{\gamma (t,t)}} 
 ( \lambda \sqrt{\gamma (t,t)} + X_t - \alpha (t) ) 
 dt
\\ 
\nonumber
 & & 
 - 
 \int_0^T 
 \int_{-\infty}^\infty 
 \phi_{\varepsilon \sqrt{\gamma (t,t)}} 
 ( y ) 
 {\bf 1}_{ 
 \{ 
 \vert 
 x - y - \alpha (t) 
 \vert 
 <  
 \lambda \sqrt{\gamma (t,t)}  
 \} }
 dy 
 dt 
\\ 
\nonumber
 & = & 
 \lambda 
 \int_{-\infty}^\infty  
 ( 
 \phi_{\varepsilon} ( - \lambda + 
 Z^{\alpha , \gamma}_t 
 ) 
 + 
 \phi_{\varepsilon} 
 ( - \lambda 
 - 
 Z^{\alpha , \gamma}_t 
 ) 
 ) 
 dt
\\ 
\nonumber
 & & 
 - 
 \int_0^T 
 \int_{-\infty}^\infty 
 \phi_{\varepsilon \sqrt{\gamma (t,t)}} 
 ( y ) 
 {\bf 1}_{\{ 
 \vert 
 x - y - \alpha (t) 
 \vert 
 <  
 \lambda \sqrt{\gamma (t,t)} 
 \} }
 dy 
 dt 
\\ 
\nonumber
 & = & 
 \lambda 
 \int_{-\infty}^\infty  
 \phi_{\varepsilon} 
 ( a - \lambda ) 
 \bar{\ell}_T^a d a  
\\ 
\nonumber
 & & 
 - 
 \int_0^T 
 \int_{-\infty}^\infty 
 \phi_{\varepsilon \sqrt{\gamma (t,t)}} 
 ( y ) 
 {\bf 1}_{\{ 
 \vert 
 x - y - \alpha (t) 
 \vert 
 <  
 \lambda \sqrt{\gamma (t,t)} 
 \} } 
 dy 
 dt 
, 
\end{eqnarray} 
 which converges in $L^2(\Omega,\P_u)$ 
 to 
$$ 
 \lambda \bar{\ell}_T^\lambda 
 - 
 \int_0^T 
 {\bf 1}_{\{ 
 \vert 
 X_t - \alpha (t) 
 \vert 
 <  
 \lambda \sqrt{\gamma (t,t)} 
 \} } 
 dt 
$$
 as $\varepsilon$ tends to zero. 
\end{proof}
\section{Numerical examples} 
\label{5} 
\noindent 
 In this section we assume that $X^u$ is a 
 centered stationary Ornstein-Uhlenbeck process 
 solution of 
$$ 
 dX^u_t= - a X^u_t dt + \sigma dB_t, \qquad t \in [0,T] 
, 
$$
 with $X^u_0\sim {\cal N} \left(0,\frac{\sigma^2}{2 a}\right)$ 
 and covariance function
$\displaystyle \gamma (s,t)=\frac{\sigma^2}{2 a} e^{-a \vert t-s\vert}$, 
 $s,t \in [0,T]$, for $\sigma, a>0$. 
 As a consequence of the following proposition 
 we can take $\Lambda = [ 0, \sqrt{2\log T}]$ 
 as parameter range when $T$ is large. 
\begin{prop}
\label{prop:C(T)}
 Assume that 
 $\Vert \alpha \Vert_{L^\infty ( [ 0,\infty ) )}  < \infty$ and 
 $\Vert u \Vert_{L^\infty ( [ 0,\infty ) )} < \infty$. 
 Then for any $r > 1$ we have 
$$ 
 \lim_{T\to \infty} 
 \P_u 
 \left( 
 \sup_{t \in [0,T]} | Z_t | \leq \sqrt{2 r \log T} 
 \right) 
 = 1. 
$$
\end{prop}
\begin{proof}
 From Theorem~1.1 of \cite{weber} 
 (see also \cite{pickands}, Theorem~2.1 of \cite{qualls}, 
 and \cite{cuzick}, page~ 488) 
 there exists a universal constants $c_1,c_2>0$ 
 such that for all $\lambda, T>0$, 
$$ 
 \P_u 
 \left( 
 \sup_{t \in [0,T]} | Z_t | > \lambda \right) 
 \leq 
 c_1 
 M(2aT,c_2/\lambda) 
 \Psi \left( \lambda \right) 
, 
$$ 
 where $\Psi (x) = \int_x^\infty e^{-y^2/2}dy/\sqrt{2\pi}$ and 
 $M(2aT,c_2/\lambda)$ is the maximal cardinal 
 of all sequences ${\cal S}$ in $[0,2aT]$ such that 
$$ 
 \Vert Z_t-Z_s\Vert_{L^2(\Omega )} = \sigma \sqrt{\frac{1-e^{-a|t-s|}}{a}} 
 > 
 \frac{c_2}{\lambda} 
, \qquad 
 s,t\in {\cal S}. 
$$ 
 Setting $\lambda = \sqrt{2 r \log T}$, $r>0$, $T>1$, 
 and using the bound 
 $\Psi ( \lambda ) \leq e^{-\lambda^2/2}/( \lambda \sqrt{2\pi})$ 
 this yields, for all $T$ large enough: 
$$ 
 \P\left( 
 \sup_{t \in [0,T]} | Z_t | \leq \sqrt{2 r \log T} 
 \right) 
 \geq 
 1 - 
 c 
 \frac{ r }{\sqrt{a}} 
 T^{1 
 - 
 r 
 }  
, 
$$ 
 which tends to $1$ as $T\to \infty$ provided $r > 1$. 
\end{proof} 
 In the next figures 
 we present some numerical simulations when 
 the signal $(X_t)_{t\in [0,T]}$ is a deterministic 
 function $(u(t))_{t\in [0,T]}$ perturbed by a centered 
 Ornstein-Uhlenbeck process, with parameters $a=0.5$, 
 $\sigma=0.05$, $T=1$. 
\\ 
 
 We represent simulated samples 
 path with the optimal thresholds 
 obtained by soft thresholding, 
 the de-noised signal after hard thresholding, and 
 the corresponding risk function 
 $(\alpha , \lambda) \mapsto \SURE_\mu (X+\xi^{\alpha , \lambda }(X) )$ 
 whose minimum gives the optimal parameter value(s). 
 The hard threshold function has not been used for estimation 
 due to increased numerical instabilities linked 
 to the simulation of the local time in \eqref{pl}. 
\subsubsection*{Simple thresholding} 
 Here we take 
 $u_t = 0.2 \times \max(0, \sin ( 3 \pi t))$, 
 $\lambda (t) = \lambda \sqrt{\gamma}$, 
 and we aim at de-noising the signal around the level $\alpha (t)=0$, 
 $t\in [0,T]$. 

\begin{figure}[!ht]
\vskip.5cm
\pgfputat{\pgfxy(7.5,-2.5)}{\pgfbox[center,center]{\pgfimage[height=60mm,width=120mm]{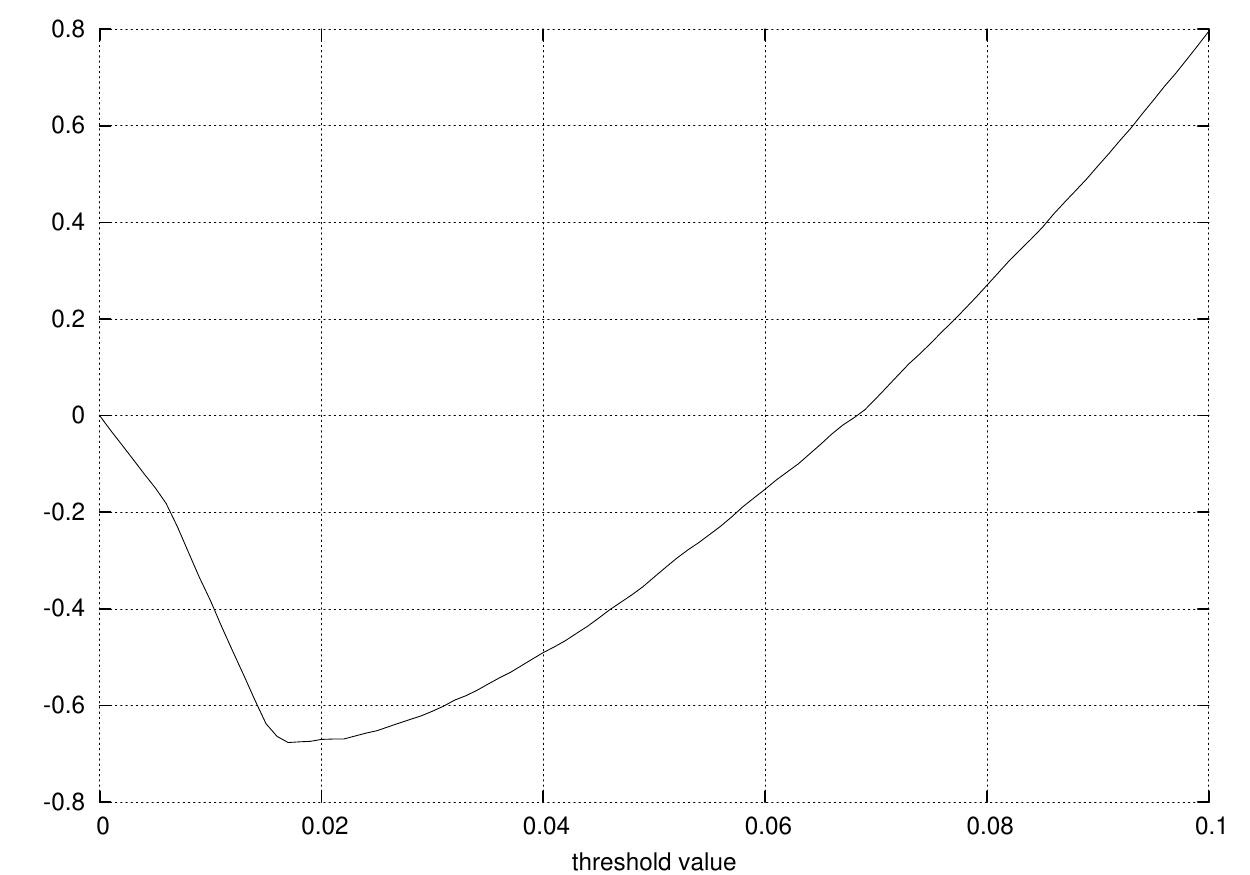}}}
\centering 
\vskip5.5cm
\caption{ 
 \small Risk function 
 $\lambda \mapsto \SURE_\mu ( X + \xi^{0 , \lambda} (X) )$ } 
\label{fs3} 
\end{figure} 

 From Figure~\ref{fs3} we estimate the 
 optimal threshold to $\lambda^* \sqrt{\gamma} = 0.018$, 
 after numerical minimization on a grid, 
 which leads to the thresholding described in Figure~\ref{fa1} 
 below. 

\begin{figure}[!ht]
\pgfputat{\pgfxy(3.5,-2.5)}{\pgfbox[center,center]{\pgfimage[height=50mm,width=74mm]{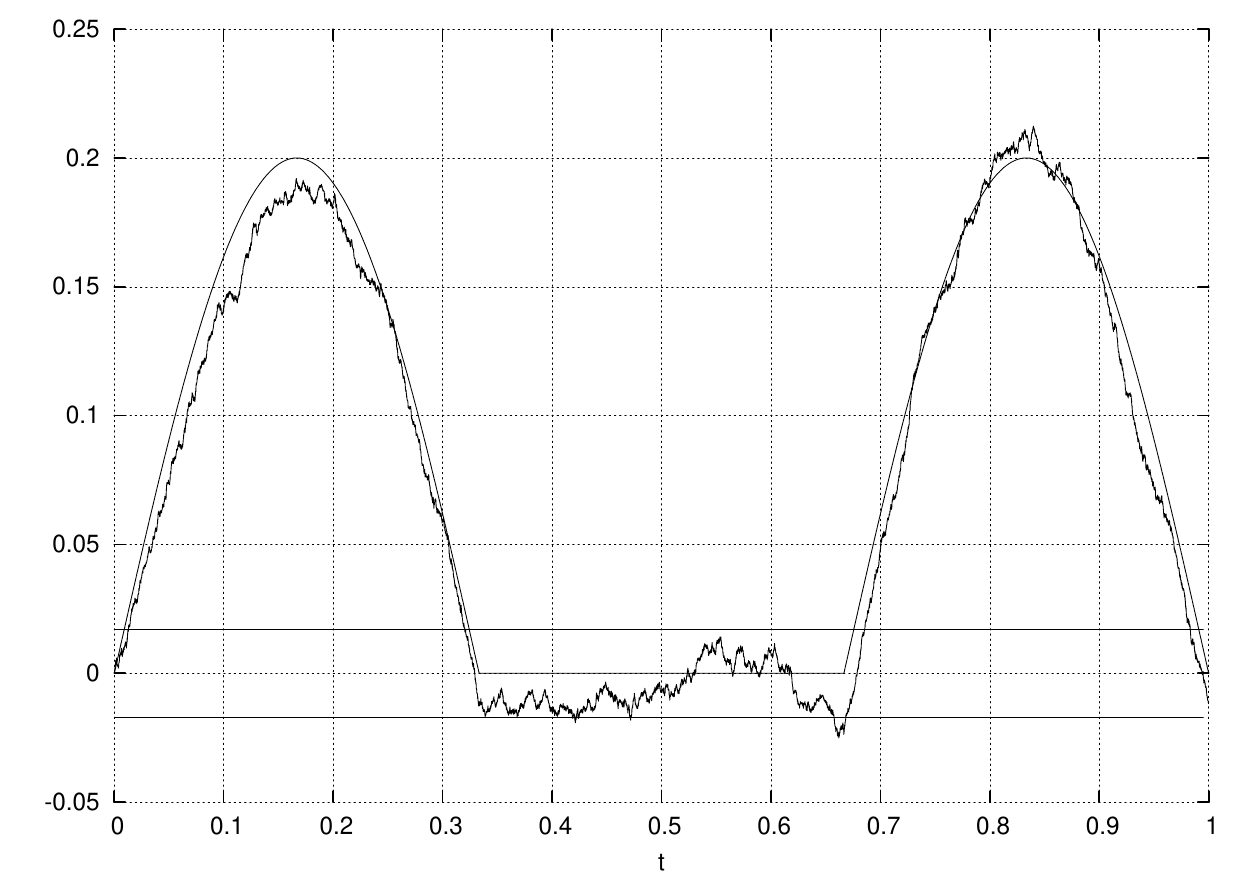}}}
\pgfputat{\pgfxy(11.6,-2)}{\pgfbox[center,center]{\pgfimage[height=50mm,width=74mm]{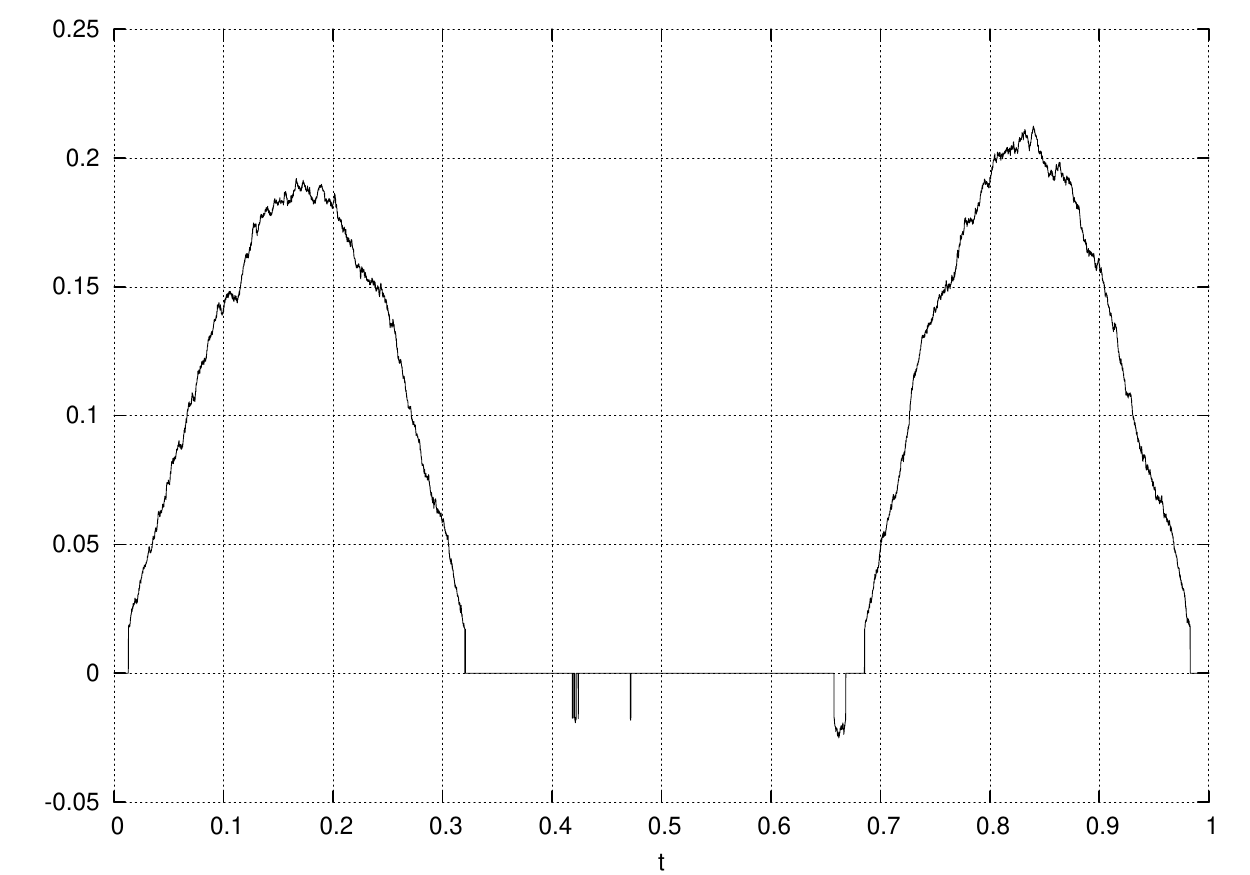}}}
\centering 
\vskip5cm
\caption{ 
 \small Process trajectory 
 ~~~~ ~~~~ ~~~~~~~~~~~~~~ 
 Estimated trajectory 
} 
\label{fa1} 
\end{figure} 

\newpage 

\subsubsection*{Level detection and thresholding} 
 We apply our method to the joint estimation 
 of parameters $\alpha$, $\lambda$, in case 
 $u_t = 0.3 + 0.2 
 \sgn (\sin (2\pi t)) \times 
 \max(0, \sin ( 3 \pi t))$, 
 $\alpha (t) = \alpha$
 and $\lambda (t) = \lambda \sqrt{\gamma}$, 
 i.e. we aim at detecting simultaneously 
 the level $\alpha = 0.3$ 
 and the threshold $\lambda \sqrt{\gamma}$ 
 at which the noise can be removed. 
 For this we have the following proposition that 
 completes Proposition~\ref{minnum}. 
\begin{prop} 
\label{minnum2} 
 The function $(\alpha , \lambda ) \mapsto 
 \SURE_\mu ( X + \xi^{\alpha , \lambda} (X) )$ 
 is continuously differentiable. 
\end{prop} 
\begin{Proof} 
 We have 
$$ 
 \frac{\partial}{\partial \alpha} 
 \SURE_\mu ( X+\xi^{\alpha , \lambda} (X) 
 ) 
 = 
 - 
 2 
 \int_0^T 
 \frac{X_t - \alpha}{\gamma (t,t)}
 {\bf 1}_{ \{ \vert X_t - \alpha \vert \leq 
 \lambda \sqrt{\gamma(t,t)} \} } 
 dt 
 + 
 2 \ell_T^{\alpha,\lambda} 
 - 
 2 \ell_T^{\alpha,-\lambda} 
, 
$$ 
 where 
 $\ell_T^{\alpha,\lambda}$ denotes the local time at level 
 $\alpha$ of the process $(X_t+\lambda \sqrt{\gamma (t,t)}))_{t\in [0,T]}$. 
\end{Proof} 

\begin{figure}[!ht]
\vskip-0.5cm
\pgfputat{\pgfxy(7.5,-2)}{\pgfbox[center,center]{\pgfimage[height=60mm,width=120mm]{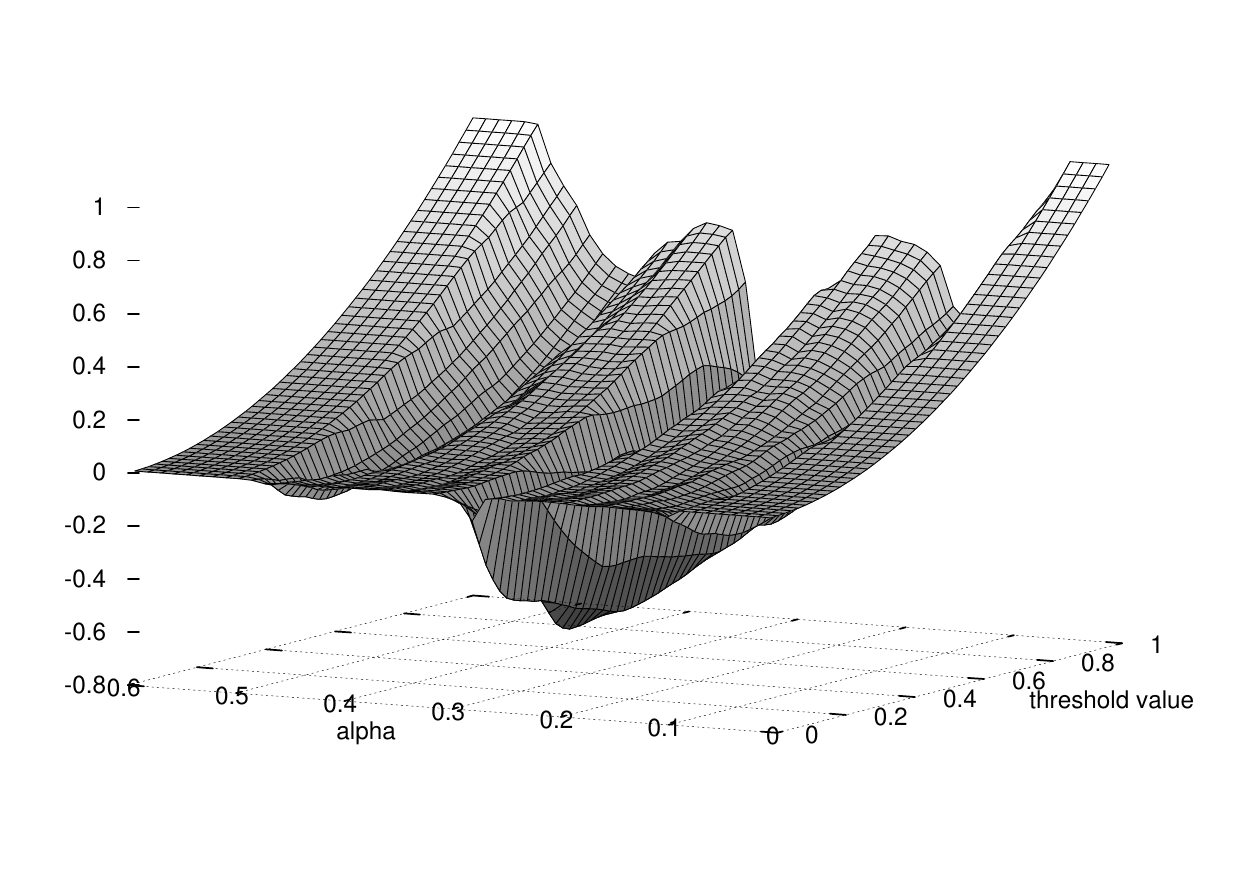}}}
\centering 
\vskip4cm
\caption{ 
 \small Risk function $(\alpha , \lambda ) \mapsto 
 \SURE_\mu ( X + \xi^{\alpha , \lambda} (X) )$ } 
\label{fa3} 
\end{figure} 
 
\newpage

 From Figure~\ref{fa3} we estimate the optimal threshold and shift parameters 
 at $\lambda^* \sqrt{\gamma} = 0.017$ and $\alpha^*=0.30$, which 
 leads to the thresholding described in Figure~\ref{fs1} below. 

\begin{figure}[!ht]
\pgfputat{\pgfxy(3.5,-2.5)}{\pgfbox[center,center]{\pgfimage[height=50mm,width=74mm]{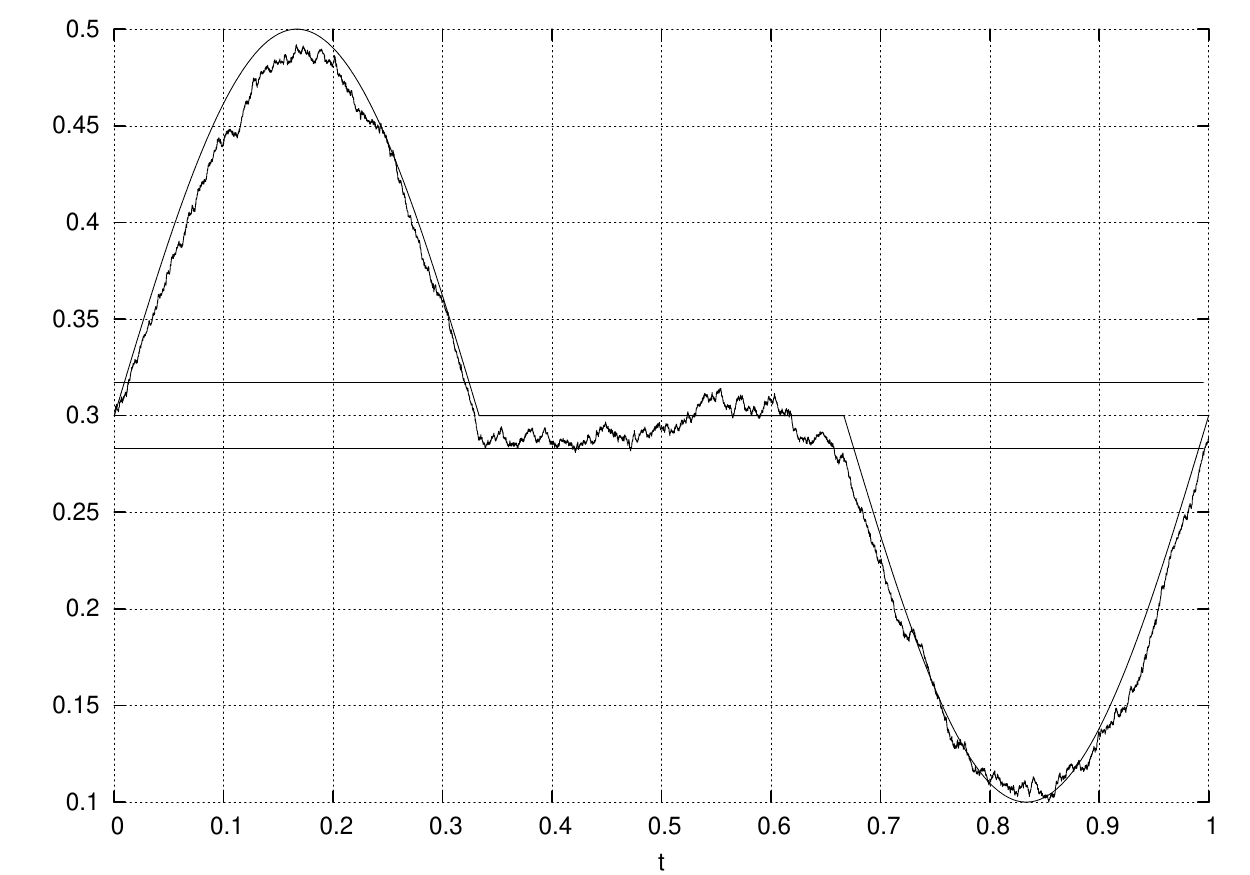}}}
\pgfputat{\pgfxy(11.6,-2)}{\pgfbox[center,center]{\pgfimage[height=50mm,width=74mm]{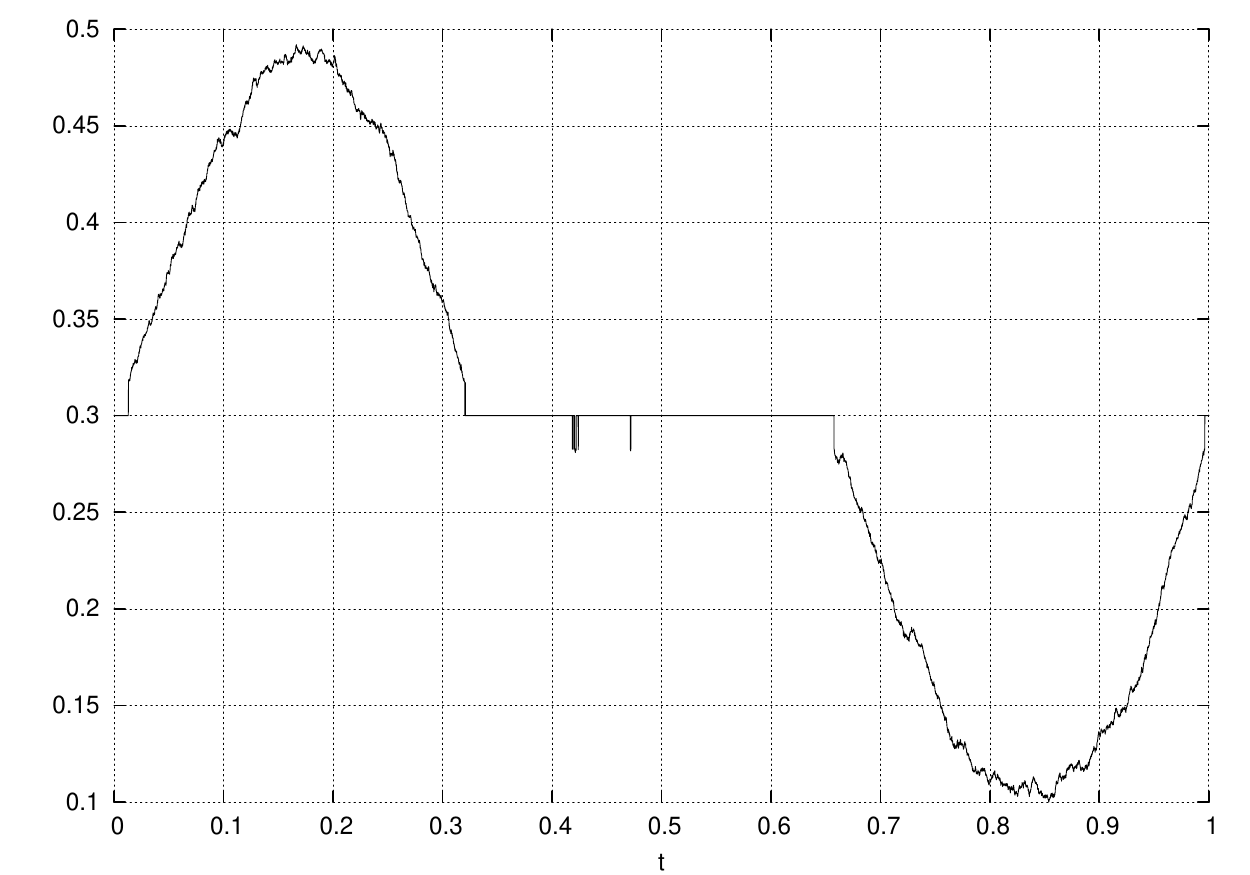}}}
\centering 
\vskip5cm
\caption{ 
 \small Process trajectory 
 ~~~~ ~~~~ ~~~~~~~~~~~~~~ 
 Estimated trajectory 
} 
\label{fs1} 
\end{figure} 

 Figure~\ref{fa3} also shows that the values $0.5$ and 
 $0.1$ are other candidates to an estimation of $\alpha$. These values 
 correspond to the extrema in the sample trajectory of Figure~\ref{fs1}. 
\subsubsection*{Drift detection and thresholding} 
 We apply our method to the joint estimation 
 of parameters $\alpha$, $\lambda$, in case 
 $u_t = 0.3 t + 0.2 \sgn (\sin (2\pi t)) 
 \times \max(0, \sin( 3 \pi t))$, 
 $\alpha (t) = \alpha t$, 
 and $\lambda (t) = \lambda \sqrt{\gamma}$, i.e. we aim at locating noise 
 with threshold $\lambda \sqrt{\gamma}$ around a line of slope $\alpha =0.3$. 
 Analogously to Propositions~\ref{minnum} and 
 \ref{minnum2} we have the following result. 
\begin{prop} 
 The function $(\alpha , \lambda ) \mapsto 
 \SURE_\mu ( X + \xi^{\alpha , \lambda} (X) )$ 
 is continuously differentiable. 
\end{prop} 
\begin{Proof} 
 We have 
$$ 
 \frac{\partial}{\partial \alpha} 
 \SURE_\mu ( X+\xi^{\alpha , \lambda} (X) 
 ) 
 = 
 - 
 2 
 \int_0^T 
 \frac{X_t - \alpha t}{\gamma (t,t)}
 {\bf 1}_{ \{ \vert X_t - \alpha t \vert \leq 
 \lambda \sqrt{\gamma(t,t)} \} } 
 t dt 
 + 
 2 \ell_T^{\alpha,\lambda} 
 - 
 2 \ell_T^{\alpha,-\lambda} 
, 
$$ 
 where $\ell_T^{\alpha,\lambda}$ denotes the local time at level  
 $\alpha$ of the process $((X_t+\lambda \sqrt{\gamma (t,t)})/t)_{t\in [0,T]}$. 
\end{Proof} 

\newpage 

\begin{figure}[!ht]
\pgfputat{\pgfxy(7.5,-2.5)}{\pgfbox[center,center]{\pgfimage[height=60mm,width=120mm]{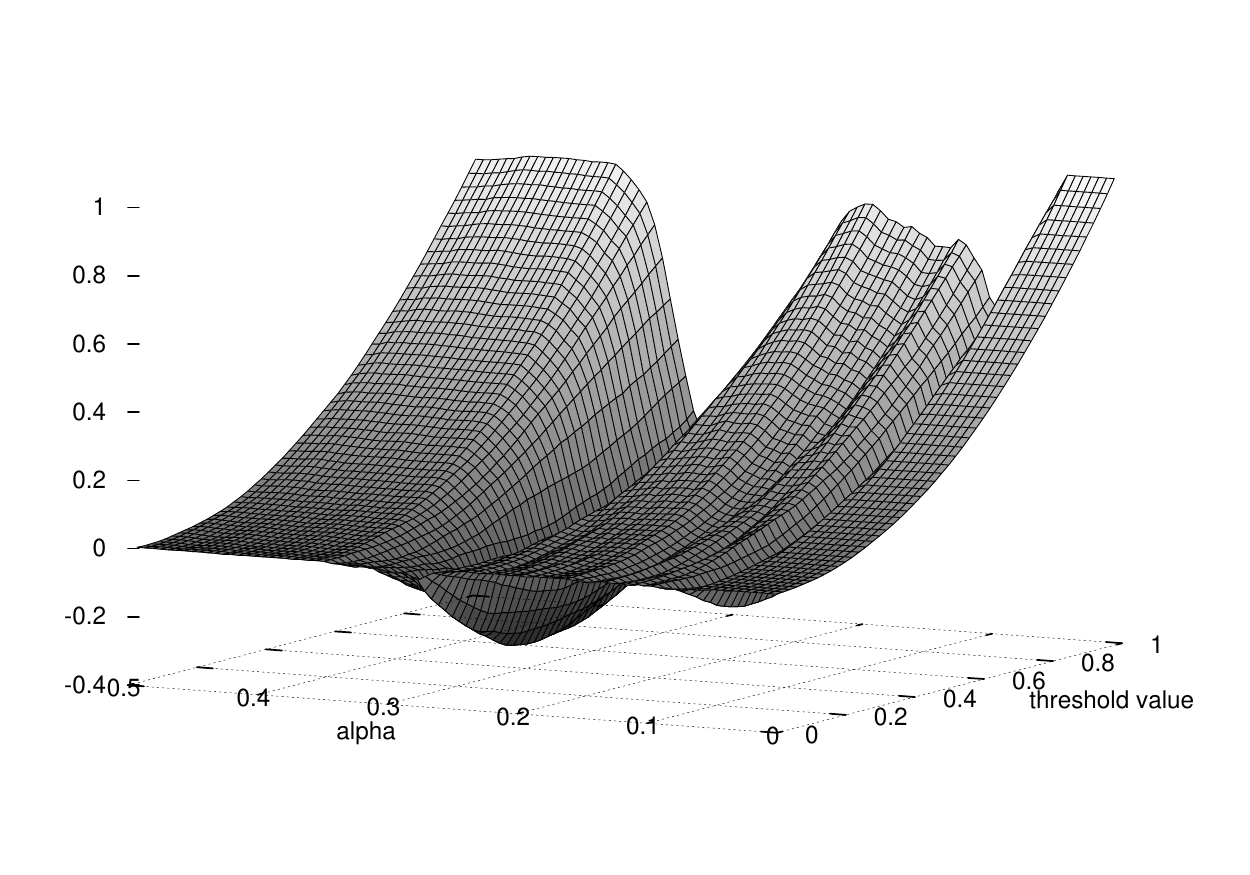}}}
\centering 
\vskip5cm
\caption{ 
 \small Risk function $(\alpha , \lambda ) \mapsto 
 \SURE_\mu ( X + \xi^{\alpha , \lambda} (X) )$ } 
\label{fl3} 
\end{figure} 
 The optimal threshold and slope parameters 
 are numerically 
 estimated at $\lambda^* \sqrt{\gamma} = 0.0093$ and $\alpha^*=0.294$. 

\begin{figure}[!ht]
\pgfputat{\pgfxy(3.5,-2.5)}{\pgfbox[center,center]{\pgfimage[height=50mm,width=74mm]{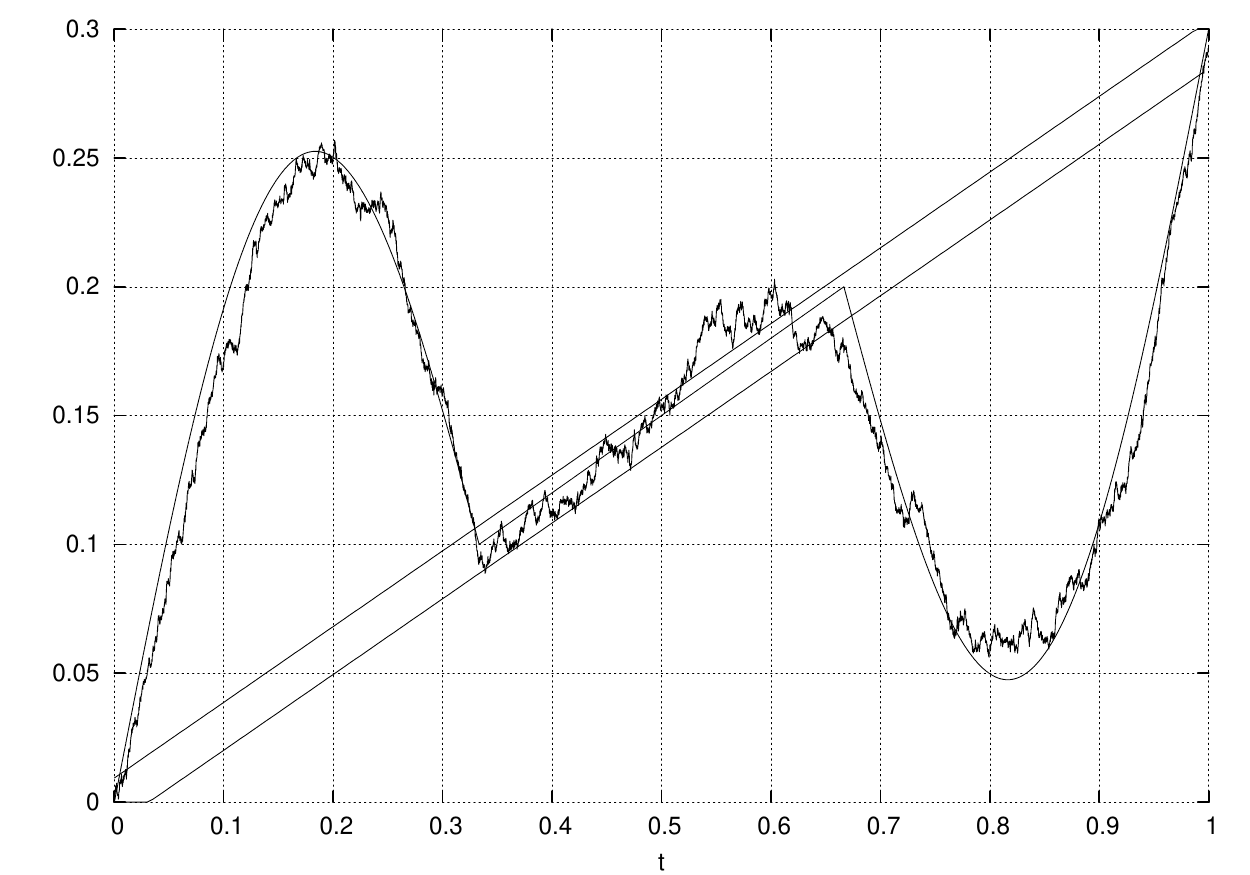}}}
\pgfputat{\pgfxy(11.6,-2)}{\pgfbox[center,center]{\pgfimage[height=50mm,width=74mm]{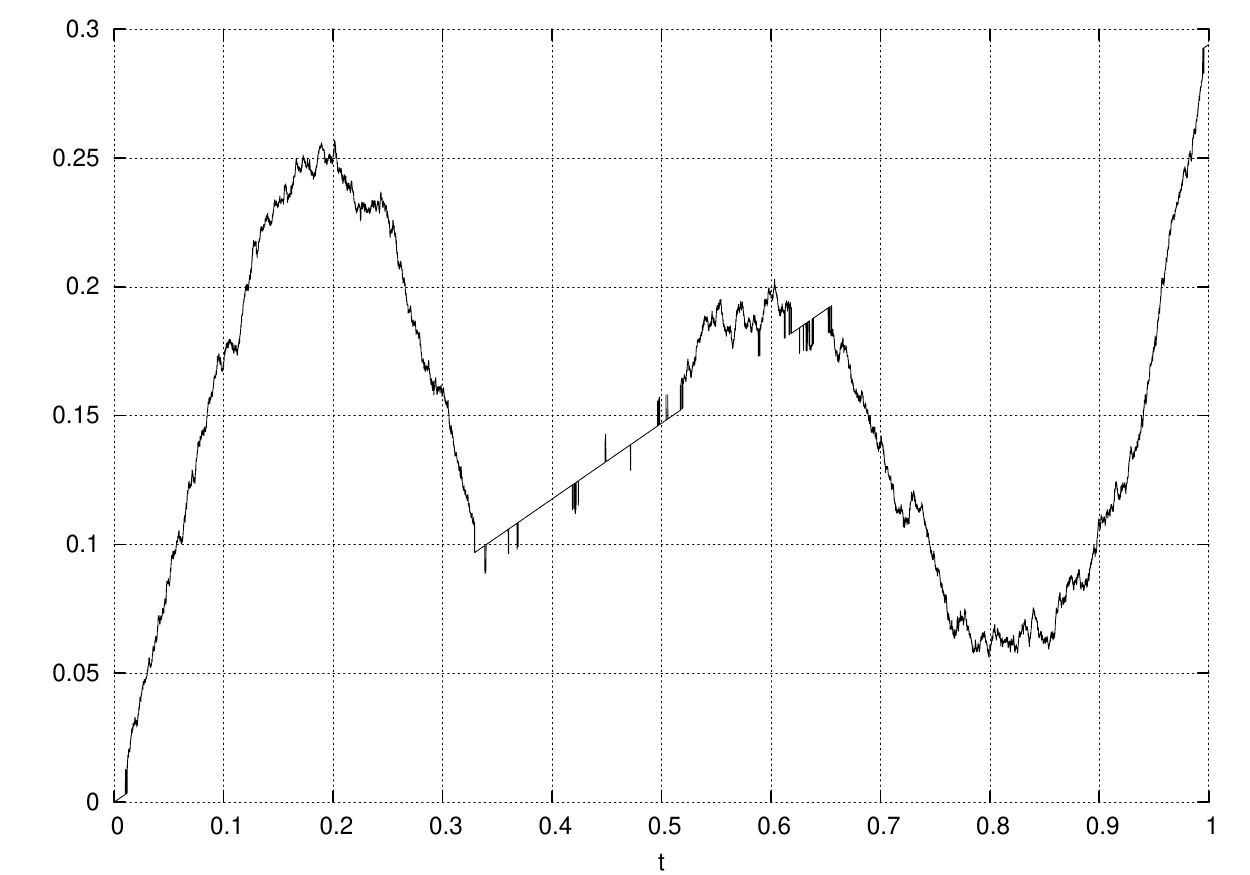}}}
\centering 
\vskip5cm
\caption{ 
 \small Process trajectory 
 ~~~~ ~~~~ ~~~~~~~~~~~~~~ 
 Estimated trajectory 
} 
\label{fl1} 
\end{figure} 

 The threshold and slope and actually slightly underestimated, 
 as the larger noise at the right end of the slope line has been interpreted as 
 being part of the signal. 

\section{Appendix} 
\label{2} 
 In this section we review three aspects of stochastic 
 analysis for Gaussian processes, including 
 local time and the Malliavin calculus calculus. 
\subsubsection*{Malliavin calculus on Gaussian space} 
 Here we recall some elements of the Malliavin calculus 
 on Gaussian space for the centered Gaussian process $(X_t)_{t\in [0,T]}$, 
 see e.g. \cite{nualartm2}. 
 Let $\mu$ be a finite Borel measure on $[0,T]$ 
 and let $\Gamma$ the operator defined as 
$$ 
 ( \Gamma g ) (t)
 =
 \int_0^T g(s) \gamma (s,t) \mu (ds ),
\qquad
 t\in [0,T]
, 
$$ 
 on the Hilbert space $H$ of functions on $[0,T]$ 
 with the inner product 
$$ 
 \langle
 h
 ,
 g
 \rangle_H
 =
 \langle
 h
 ,
 \Gamma
 g
 \rangle_{L^2([0,T],d\mu )}
. 
$$

 The process $(X_t)_{t\in [0,T]}$ can be used to construct 
 an isometry $X : H \to L^2 (\Omega , {\cal F} , P)$ as 
$$ 
 X (h) = \int_0^T X_s h(s) \mu (ds), \qquad h\in H
. 
$$ 
 Then $\{ X (h) \ : \ h\in H\}$ is an 
 isonormal Gaussian process on $H$, i.e. 
 a family of centered
 Gaussian random variables satisfying
$$\E
 [ X (h) X (g) ]
 = \langle h , g \rangle_H,
 \qquad
 h,g \in H
. 
$$
 For any orthonormal basis $(h_k)_{k\in\inte}$ of
 $L^2([0,T],d\mu)$, we have the Karhunen-Lo\`eve expansion 
\begin{equation} 
\label{kl} 
 X_t = \sum_{k=0}^\infty h_k (t) X (h_k)
,
 \qquad
 t\in [0,T]
. 
\end{equation} 

 Let now ${\cal S}$ denote the space of cylindrical functionals of the
 form 
\begin{equation}
\label{e1}
 F= f_n
 \left(
 X^u ( h_1 )
 ,
 \ldots
 ,
 X^u ( h_n)
 \right)
,
\end{equation}
 where $f_n$ is in the space of infinitely differentiable rapidly decreasing
 functions on $\real^n$, $n\geq 1$.
\noindent
\begin{definition}
 The $H $-valued Malliavin derivative is defined as
$$
 \nabla_t F =
 \sum_{i=1}^{n}
 h_i ( t )
 \partial_i f_n
 \left(
 X^u (h_1 )
 ,
 \ldots
 ,
 X^u (h_n)
 \right)
,
$$
 for $F \in {\cal S}$ of the form \eqref{e1}.
\end{definition}
\noindent
 It is known that $\nabla$ is closable, cf. Proposition~1.2.1 of
 \cite{nualartm2}, and its closed domain will be denoted by $\Dom (\nabla )$.
\begin{definition} 
\label{defd} 
 Let $D_t$ be defined on $F\in \Dom (\nabla )$ as
$$D_t F : = ( \Gamma \nabla F ) (t)
, \quad t\in [0,T].
$$ 
\end{definition}
\noindent
 Let $\delta : L^2_u (\Omega ; H ) \to L^2 (\Omega , \P_u )$
 denote the closable adjoint of $\nabla$, i.e. the
 divergence operator under $\P_u$,
 which satisfies the integration by parts formula
\begin{equation}
\label{intparpartiesvar1/N}
 \E_u [ F \delta ( v ) ]
 = \E_u [ \langle v , \nabla F \rangle_{H} ]
,
 \qquad
 F \in \Dom ( \nabla )
 ,
 \quad
 v \in \Dom (\delta )
, 
\end{equation} 
 where $\E_u$ denotes the expectation under $\P_u$, 
 with the relation
$$ 
 \delta (h F ) = F X (h) - \langle h , \nabla F\rangle_H,
$$
 cf. \cite{nualartm2}, for $F\in \Dom (\nabla )$ and $h\in H$ such
 that $hF\in \Dom (\delta )$. 
 The next lemma will be needed 
 in Proposition~\ref{lemma1} below 
 to establish 
 Stein's Unbiased Risk Estimate for Gaussian 
 processes.
\begin{lemma}
\label{prt} 
 For any $F \in \Dom ( \nabla )$ and $u\in H$ we have
$$ 
 \E_u [ F X^u_t ]
 = \E_u [
 D_t F
 ]
,
 \qquad
 t\in [0,T]
.
$$
\end{lemma}
\begin{Proof}
 We have
\begin{eqnarray*}
 \E_u [ F X^u_t ]
 & = &
 \sum_{k=0}^\infty
 h_k(t)
 \E_u [ F X^u (h_k) ]
\\
 & = &
 \sum_{k=0}^\infty
 h_k(t)
 \E_u [ F \delta (h_k) ]
\\
 & = &
 \sum_{k=0}^\infty
 h_k(t)
 \E_u [ \langle h_k , \nabla F\rangle_H ]
\\
 & = &
 \sum_{k=0}^\infty
 h_k(t)
 \E_u [ \langle h_k , \Gamma \nabla F\rangle_{L^2([0,T],\mu )} ]
\\
 & = &
 \E_u [
 ( \Gamma \nabla F ) (t)
 ]
,
 \qquad
 F \in \Dom ( \nabla )
 ,
 \quad
 t\in [0,T]
.
\end{eqnarray*}
\end{Proof}
\noindent
 Note that since $u\in H$ we have 
 $\nabla_s X_t (h) = \nabla_s X^u_t (h) = h(s)$ and 
\begin{eqnarray*} 
 D_t X_t & = & 
 ( \Gamma \nabla X_t ) (t)
\\ 
 & = & 
\int_0^T \gamma(s,t) \nabla_s X_t \mu(ds) 
\\ 
 & = & 
\sum_{k=0}^\infty h_k(t) \int_0^T \gamma(s,t) \nabla_s X(h_k) \mu(ds) 
\\ 
 & = & 
\sum_{k=0}^\infty h_k(t) \langle \gamma(\cdot,t),h_k\rangle_{L^2([0,T],d\mu)} 
\\ 
 & = & 
 \gamma(t,t), 
 \quad t\in [0,T]
.
\end{eqnarray*} 

\subsubsection*{Local time of Gaussian processes} 
 Given $(Z_t)_{t\in [0,T]}$ a Gaussian process let 
$$ 
 \Delta (s,t) = \var ( Z_t - Z_s), 
 \qquad 
 0\leq s , t \leq T 
, 
$$ 
 and denote by 
$$ 
 L_T^\lambda 
 : = 
 \int_0^T 
 {\bf 1}_{ \{ 
 Z_t \leq \lambda 
 \} 
 } dt 
$$ 
 the occupation time of $(Z_t)_{t\in [0,T]}$ 
 up to $T$ in the set $( - \infty , \lambda ]$. 
\\ 
 
 Recall that a classical result of Berman~\cite{berman}, see 
 Theorem~21.9 of \cite{GemanHorowitz}, shows that if 
\begin{equation} 
\label{eq:TL} 
\int_0^T \int_0^T 
 \Delta^{-1} (s,t) ds dt<\infty
, 
\end{equation} 
 then for any $\lambda\in \real$ the local time 
$$ 
 \ell_T^\lambda : = 
 \frac{\partial}{\partial \lambda} L_T^\lambda 
$$ 
 of $(Z_t)_{t\in [0,T]}$ at the level $\lambda$ exists 
 and the occupation time density formula 
\begin{equation} 
\label{otdf} 
 \int_0^T f ( Z_t) dt 
 = 
 \int_\real f ( \lambda ) \ell_T^\lambda d \lambda 
\end{equation} 
 holds for every positive measurable function $f$ on 
 $\real$. 
 The local time $\bar{\ell}_T^\lambda$ 
 of $|Z_t|$ is given by 
 $\bar{\ell}_T^\lambda = 
 \ell_T^{-a} 
 + 
 \ell_T^a$ 
 and the related occupation time 
 formula can be obtained under 
 the same condition from the relation 
$$ 
 \int_0^T f ( 
 | Z_t
 | 
 ) 
 dt 
 = 
 \int_{-\infty}^\infty 
 f(|a|) 
 \ell_T^a 
 da 
 =  \int_0^\infty 
 f(a) 
 \bar{\ell}_T^a 
 da 
. 
$$

\end{document}